%% file: TrajFlat.tex
\documentclass[10pt,conference]{ieeeconf}

\usepackage{amsmath,amssymb,graphicx}
\usepackage{url}
\usepackage{algorithm}
\usepackage{algpseudocode}
\usepackage{booktabs}

\author{Basile Graf, Gustave Lapierre, and Philippe M\"ullhaupt}

\title{Efficient flatness-based computation of trajectories for vehicles with many trailers}

\begin{document}
\maketitle

\begin{abstract}
We address the practical computation of feedforward steering inputs for a car with $n$ trailers via differential flatness. The standard route---iterated symbolic differentiation of the flat output---scales prohibitively with~$n$: even with computer-algebra software the expression sizes grow super-exponentially and become intractable beyond $n\approx 5$. We propose an algorithm that combines four ingredients (angular intermediate variables, factorisation of the recursion in sub-maps, rescaled higher-order product/quotient rules, and a composition rule based on truncated formal power series instead of Fa\`a di Bruno's formula). We benchmark the proposed method against (i) direct symbolic differentiation in SymPy and (ii) Fa\`a di Bruno's formula evaluated via Bell polynomials. Both reference methods hit a clear computational wall: SymPy reaches a $20\,\mathrm{s}$ timeout at derivative order $r\!=\!14$, while the Bell-polynomial recursion exceeds a $60\,\mathrm{s}$ timeout beyond $r\!=\!24$. The proposed approach has $O(r^{3})$ arithmetic complexity in the derivative order and remains under one millisecond up to $r\!=\!40$, with numerical agreement to machine precision in the regime where the reference methods succeed. As an end-to-end illustration, a $20$-trailer parking manoeuvre is computed and animated.
\end{abstract}

\section{Introduction}

Differential algebra as presented in \cite{Ritt} and \cite{Kolchin} was recognised in \cite{Fliess1} as providing valuable insight for control systems, and the resulting definitions of control properties can be extended to a class of nonlinear systems \cite{Fliess2}. The notion of \emph{differential flatness} was first introduced using the language of differential algebra in \cite{FLMR1}. Shortly after, a differential-geometric description in terms of diffieties was given in \cite{FLMR2,FLMR3}. Necessary and sufficient conditions in the nonlinear geometric setting were obtained in \cite{Levine1}, and an elegant treatment for linear systems in~\cite{LevineNguyen}.

For our purposes we adopt the following informal working definition. Consider a system in explicit form
\begin{equation}
\label{eq:syst}
\dot x_i = f_i(x_1,\ldots,x_n,u_1,\ldots,u_m), \qquad i=1,\ldots,n.
\end{equation}
The system is \emph{flat} if there exist $m$ smooth functions
\[
z_s = z_s\bigl(t,x,u,\dot u,\ldots,u^{(L-1)}\bigr),\qquad s=1,\ldots,m,
\]
called \emph{flat outputs}, for some finite~$L$, together with a map $\Phi:\mathbb{R}^{1+mR}\to\mathbb{R}^{n}$ such that
\[
x_i = \Phi_i\bigl(t,z,\dot z,\ldots,z^{(R-1)}\bigr),\qquad i=1,\ldots,n,
\]
for some finite~$R$, with an analogous expression for the inputs. Trajectory planning in flat coordinates then reduces to designing a smooth curve $t\mapsto z(t)$; the corresponding state and input trajectories are recovered by symbolic substitution and differentiation. The bottleneck is exactly this differentiation step when many derivatives are required.

\paragraph*{Contributions and outline.} For the canonical car-with-trailers benchmark of \cite{FLMR0}, naive symbolic differentiation in computer-algebra packages such as SymPy, Maple or Mathematica becomes the dominant cost beyond a handful of trailers. Fa\`a di Bruno's formula evaluated through Bell polynomials \cite{Faa,wiki1,wiki2} is sometimes proposed as a way out, but the partition-based recursion is itself super-polynomial. We document both walls quantitatively and present an alternative that we have used to compute flat-output derivatives for up to $n=20$ trailers. Section~II recalls the flat parametrisation of the trailer chain. Section~III details the proposed four-step algorithm. Section~IV reports the quantitative comparison between direct symbolic differentiation, Bell-polynomial composition, and the proposed formal-series composition. Section~V illustrates the end-to-end pipeline with simulations. Section~VI concludes.

Throughout, $'$ denotes the derivative with respect to the path parameter~$s$. Bold lower-case symbols denote two-dimensional column vectors. By $\arctan_2(y,x)$ we mean the standard two-argument arc-tangent that returns the angle of $(x,y)\in\mathbb{R}^2\setminus\{0\}$ in $(-\pi,\pi]$.

\section{Differential flatness of a vehicle with trailers}

Explicit computation of the steering input for non-holonomic vehicles is a classical application of differential flatness. In particular, it allows one to determine the required steering profile of a car with $n$ trailers given the path profile that the last trailer's axle should follow.

\begin{figure}[h]
\begin{center}
\includegraphics[width=0.55\textwidth]{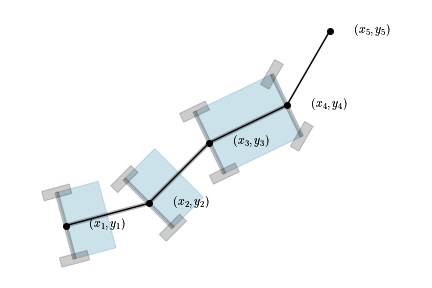}
\end{center}
\caption{A car with two trailers. The positions of the middle of the axle of each trailer are $(x_1,y_1)$, $(x_2,y_2)$, $(x_3,y_3)$ and $(x_4,y_4)$. The numbering starts with the trailing cart (number~1). The point $(x_5,y_5)$ is fictitious (not physical) and is introduced solely to encode the orientation of the front car's steering wheels.}
\label{fig:train}
\end{figure}

For instance, consider a car with two trailers (Fig.~\ref{fig:train}). One may choose any sufficiently differentiable path $\bigl(x_1(s),y_1(s)\bigr)^{T}$ for the axle of the last trailer. From this path and its derivatives, all the other axle paths $\bigl(x_k(s),y_k(s)\bigr)^{T}$, $k=2,\ldots,n$, can be computed explicitly. The angle between the link~$n{-}1$\,--\,$n$ and the link~$n$\,--\,$n{+}1$ is the required steering angle of the head car. Note also that the four-wheeled head car is equivalent to a two-wheeled trailer, which is why we use a single set of formulas throughout.

The equations are elementary, but naively applying them directly leads to unwieldy expressions whose size explodes with~$n$ (Section~IV).

\subsection{Recursive kinematic equation}

The kinematic constraint that the next trailer follows the heading of the current one \cite{FLMR0} yields, in path-parameter form,
\begin{equation}
\label{eq:kinematic}
\begin{pmatrix} x_{k+1}\\ y_{k+1}\end{pmatrix}
= \begin{pmatrix} x_{k}\\ y_{k}\end{pmatrix}
+ \frac{L_k}{\sqrt{\bigl(x_k'(s)\bigr)^{2} + \bigl(y_k'(s)\bigr)^{2}}}
\begin{pmatrix} x_k'(s)\\ y_k'(s)\end{pmatrix},
\end{equation}
where $L_k>0$ is the inter-axle distance between trailer~$k$ and trailer~$k+1$. Defining $\boldsymbol{q}_k(s)\triangleq\bigl(x_k(s),y_k(s)\bigr)^{T}$ and iterating \eqref{eq:kinematic} $n-1$ times, one obtains a map
\begin{equation}
\label{eq:phi}
\varphi:\;\bigl(\boldsymbol{q}_1(s),\boldsymbol{q}_1'(s),\ldots,\boldsymbol{q}_1^{(n-1)}(s)\bigr)
\;\longmapsto\;
\bigl(\boldsymbol{q}_1(s),\ldots,\boldsymbol{q}_n(s)\bigr)
\end{equation}
that returns all trailer positions given the derivatives of the chosen flat output $\boldsymbol{q}_1(s)$. Applying~\eqref{eq:kinematic} a few times in a computer-algebra system already produces very large expressions, mainly because of the repeated differentiation of the square-root term.

\subsection{Angle variables}

The first ingredient is to introduce angle variables and rewrite the recursion as
\begin{align}
\label{eq:angle-recursion}
\begin{pmatrix} x_{k+1}\\ y_{k+1}\end{pmatrix} &=
\begin{pmatrix} x_{k}\\ y_{k}\end{pmatrix}
+ L_k\begin{pmatrix} \cos\alpha_k\\ \sin\alpha_k\end{pmatrix},\\
\label{eq:angle-def}
\alpha_k &\triangleq \arctan_2\bigl(y_k'(s),\,x_k'(s)\bigr).
\end{align}
Differentiation of \eqref{eq:angle-def} with respect to $s$ yields the \emph{rational} expression
\begin{equation}
\label{eq:alpha-prime}
\alpha_k'(s) = \frac{x_k'(s)\,y_k''(s)-x_k''(s)\,y_k'(s)}
{\bigl(x_k'(s)\bigr)^{2}+\bigl(y_k'(s)\bigr)^{2}},
\end{equation}
free of the square-root singularity. Higher derivatives $\alpha_k^{(r)}(s)$ are then obtained from the derivatives of $x_k(s)$ and $y_k(s)$ by repeated application of the product, quotient and chain rules. The derivatives of $\cos\alpha_k$ and $\sin\alpha_k$ are obtained by chain rule from those of $\alpha_k$. Higher-order forms of the product, quotient and chain (i.e.\ composition) rules are detailed in Sections~III.B and~III.C below.

\subsection{Factorisation of the iterated map}

The second ingredient is to factorise $\varphi$ into single-trailer steps instead of computing it as a single $n-1$-fold composition. Differentiating \eqref{eq:angle-recursion} $n-k-1$ times gives a map
\begin{eqnarray}
\label{eq:phi_k}
\phi_k:\;
\bigl(\boldsymbol{q}_k(s),\boldsymbol{q}_k'(s),\ldots,\boldsymbol{q}_k^{(n-k)}(s)\bigr)
\longmapsto \nonumber \\
\bigl(\boldsymbol{q}_{k+1},\boldsymbol{q}_{k+1}',\ldots,\boldsymbol{q}_{k+1}^{(n-k-1)}\bigr)
\end{eqnarray}
so that the full map admits the factorisation
\begin{equation}
\label{eq:phi-factorised}
\varphi = \phi_{n-1}\circ\phi_{n-2}\circ\cdots\circ\phi_{1}.
\end{equation}
Each $\phi_k$ requires one fewer derivative than the previous one and operates on \emph{numerical} arrays of derivatives, not on symbolic expressions. We pay $O(n)$ such applications, each of polynomial cost (see Section~IV).

\section{Higher-order calculus rules}

We collect here the higher-order versions of the product, quotient and composition rules used by each~$\phi_k$.

\subsection{Higher-order product, reciprocal and quotient}

For two scalar functions $f,g$ of~$s$, the standard Leibniz rule reads
\begin{equation}
\label{eq:leibniz}
\frac{d^{r}}{ds^{r}}\bigl(f\,g\bigr)
= \sum_{j=0}^{r}\binom{r}{j}\,f^{(r-j)}\,g^{(j)}.
\end{equation}
For the reciprocal $\bar g\triangleq 1/g$, we use the relation $g\cdot\bar g\equiv 1$, take the $r$-th derivative with~\eqref{eq:leibniz}, and solve for $\bar g^{(r)}$, which yields a recurrence
\begin{equation}
\label{eq:recip}
\bar g^{(r)} = -\,\frac{1}{g}\sum_{j=1}^{r}\binom{r}{j}\,g^{(j)}\,\bar g^{(r-j)}.
\end{equation}
The quotient rule follows by combining \eqref{eq:leibniz}--\eqref{eq:recip}: $\bigl(f/g\bigr)^{(r)} = \bigl(f\cdot\bar g\bigr)^{(r)}$.

\paragraph*{Rescaled derivatives.} Define
\begin{equation}
\label{eq:rescaled}
f^{[k]}(s)\triangleq\frac{1}{k!}\,f^{(k)}(s).
\end{equation}
The Leibniz rule then takes the convolutional form
\begin{equation}
\label{eq:leibniz-rescaled}
(f\,g)^{[r]} = \sum_{j=0}^{r} f^{[r-j]}\,g^{[j]},
\end{equation}
i.e.\ the coefficient-wise convolution of the truncated Taylor coefficients. Equation~\eqref{eq:leibniz-rescaled} has two advantages: (i)~the binomial coefficients of \eqref{eq:leibniz} disappear, simplifying both code and analysis; (ii)~the high-order derivatives $f^{(k)}$ grow factorially in $k$, whereas the rescaled $f^{[k]}$ are bounded for analytic~$f$, leading to much better numerical conditioning. We use \eqref{eq:leibniz-rescaled} throughout.

\subsection{Higher-order composition: two routes}

To complete the computation of each $\phi_k$, we need the higher-order derivatives of compositions of the form $\cos\circ\,\alpha_k(s)$ and $\sin\circ\,\alpha_k(s)$.

\paragraph*{Fa\`a di Bruno's formula.}
The classical formula \cite{Faa,wiki1} expresses
\begin{equation}
\label{eq:faa}
\bigl(f\circ g\bigr)^{(n)} = \sum_{k=1}^{n} f^{(k)}(g)\;B_{n,k}\bigl(g',g'',\ldots,g^{(n-k+1)}\bigr),
\end{equation}
where $B_{n,k}$ are the partial (exponential) Bell polynomials~\cite{wiki2}. A standard recurrence \cite{wiki2}
\begin{equation}
\label{eq:bell}
B_{n,k}(x_1,\ldots,x_{n-k+1}) = \sum_{i=0}^{n-k}\binom{n-1}{i}\,x_{i+1}\,B_{n-i-1,k-1}
\end{equation}
makes \eqref{eq:faa} look attractive, but a direct recursive evaluation of \eqref{eq:bell} performs a number of operations that is proportional to the number of integer partitions of~$n$, which by Hardy--Ramanujan grows as
\(
p(n)\sim \tfrac{1}{4n\sqrt{3}}\exp\!\bigl(\pi\sqrt{2n/3}\bigr).
\)
Memoisation reduces this to $O(n^{3})$ but with a large constant; the implementation we use here (see~\cite{Code}) is unmemoised, which makes the wall visible at moderate~$n$.

\paragraph*{Formal-series composition.}
A second route, advocated by several authors of the combinatorics literature \cite{combinatorics}, is to view both $f$ and $g$ as truncated formal power series at a chosen base point $s_0$, and to compose them \emph{as polynomials}. Writing the Taylor expansion of $f(g(s))$ about~$g(s_0)$,
\[
f\bigl(g(s)\bigr) = \sum_{k\ge 0} \tfrac{1}{k!}\,f^{(k)}\!\bigl(g(s_0)\bigr)\,\bigl(g(s)-g(s_0)\bigr)^{k},
\]
and substituting the Taylor series of $g$,
\[
g(s) = \sum_{l\ge 0} \tfrac{1}{l!}\,g^{(l)}(s_0)\,(s-s_0)^{l},
\]
yields $f(g(s)) = (p_1\circ p_2\circ p_3)(s-s_0)$ with
\begin{align}
p_1(x) &\triangleq \sum_k \tfrac{1}{k!}\,f^{(k)}\!\bigl(g(s_0)\bigr)\,x^{k}, \nonumber\\
p_2(y) &\triangleq -g(s_0)+y, \label{eq:p123}\\
p_3(z) &\triangleq \sum_l \tfrac{1}{l!}\,g^{(l)}(s_0)\,z^{l}. \nonumber
\end{align}
Truncating $p_1,p_2,p_3$ at degree~$r$ and composing them \emph{as polynomials} returns the first $r$ coefficients of $f\circ g$ exactly. In rescaled form $f^{[k]},g^{[l]}$, all factorial factors disappear and one manipulates polynomials with rational---in practice floating-point---coefficients only.

\paragraph*{Implementation.} The polynomial composition $p\circ q$ truncated at degree $r$ is built monomial-by-monomial,
\[
(p\circ q)\bmod s^{r+1} \;=\; \sum_{k=0}^{r} p_k\,\bigl(q^{k}\bmod s^{r+1}\bigr),
\]
where each $q^{k+1}$ is obtained from $q^{k}$ by one convolutional product~\eqref{eq:leibniz-rescaled}. This costs $O(r^{2})$ for one product and $O(r^{3})$ for the full composition, with no symbolic blow-up; the constant is small enough to make $r\!=\!40$ run in under a millisecond on a laptop (Section~IV).

\subsection{Algorithm summary}

The complete pipeline is summarised in Algorithm~\ref{alg:main}.

\begin{algorithm}[t]
\caption{Flat-output derivatives for $n$ trailers}\label{alg:main}
\begin{algorithmic}[1]
\Require Path $\boldsymbol{q}_1(s)=(x_1,y_1)$; inter-axle lengths $L_1,\ldots,L_{n-1}$; orders $r_k\!=\!n{-}k$
\State Evaluate Taylor coefficients $\boldsymbol{q}_1^{[0]},\ldots,\boldsymbol{q}_1^{[r_1]}$ at $s_0$
\For{$k=1,\ldots,n-1$}
   \State Compute $(\alpha_k')^{[0]},\ldots,(\alpha_k')^{[r_k-1]}$ from \eqref{eq:alpha-prime} \par
   \hskip\algorithmicindent using \eqref{eq:leibniz-rescaled} for products and \eqref{eq:recip} for the reciprocal
   \State Integrate one term to get $\alpha_k^{[0]},\ldots,\alpha_k^{[r_k]}$
   \State Compute $(\cos\alpha_k)^{[0]},\ldots,(\cos\alpha_k)^{[r_k]}$ and $(\sin\alpha_k)^{[\bullet]}$ \par
   \hskip\algorithmicindent by formal-series composition of \eqref{eq:p123}
   \State $\boldsymbol{q}_{k+1}^{[\bullet]} \gets \boldsymbol{q}_{k}^{[\bullet]} + L_k\,\bigl((\cos\alpha_k)^{[\bullet]},(\sin\alpha_k)^{[\bullet]}\bigr)$
\EndFor
\Ensure $\boldsymbol{q}_2^{[\bullet]},\ldots,\boldsymbol{q}_n^{[\bullet]}$ at $s_0$
\end{algorithmic}
\end{algorithm}

\section{Complexity comparison and numerical validation}
\label{sec:bench}

We compare three implementations of the elementary block underlying~\eqref{eq:phi_k}: the computation of the first~$r$ derivatives of $h(s)=\cos\bigl(g(s)\bigr)$ at a fixed point $s_0\in\mathbb{R}$, where $g$ is a polynomial of degree~$r$ with integer coefficients (drawn once for each~$r$ with a fixed pseudo-random seed). The three methods are:
\begin{description}
  \item[(M1)] \emph{Direct symbolic differentiation.} Repeated calls to \texttt{sympy.diff} on $\cos(g(s))$ without simplification, using SymPy~1.14.
  \item[(M2)] \emph{Fa\`a di Bruno via Bell polynomials.} The \texttt{composition\_rule} of \cite{Code} that implements \eqref{eq:faa}--\eqref{eq:bell}.
  \item[(M3)] \emph{Formal-series composition (proposed).} The \texttt{composition\_rule\_r} of \cite{Code} that implements \eqref{eq:p123}.
\end{description}
Choosing $r$ as the figure of merit is natural because the $k$-th trailer in the chain~\eqref{eq:phi-factorised} needs derivatives up to order~$r_k=n-k$, cf.~\eqref{eq:phi_k}, so the per-trailer cost is dominated by the case $r\approx n$; the benchmark range $r\le 40$ therefore covers twice the largest order arising in the $20$-trailer manoeuvre of Section~V.

\paragraph*{Methodology.} All experiments were run on a single core of a 2024 laptop using Python~3.12 with SymPy~1.14 and NumPy~2.2. Per-row timeouts were set to $20\,\mathrm{s}$ for~(M1) and $60\,\mathrm{s}$ for~(M2). For~(M1) we report the number of operations \texttt{count\_ops} of the highest-order derivative produced as an order-of-magnitude proxy for the expression size; we skip the substitution $s\mapsto s_0$ once the expression exceeds $10^{5}$ nodes (the differentiation cost is unaffected).

\begin{figure}[t]
\centering
\includegraphics[width=0.95\columnwidth]{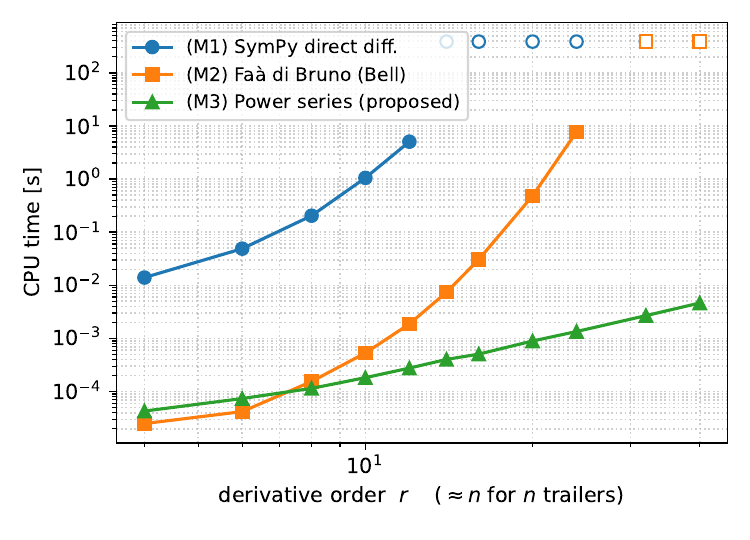}
\caption{CPU time of the three methods as a function of the derivative order~$r$ (log--log). Open markers indicate timeouts (20\,s for M1, 60\,s for M2).}
\label{fig:bench-time}
\end{figure}

\begin{figure}[t]
\centering
\includegraphics[width=0.95\columnwidth]{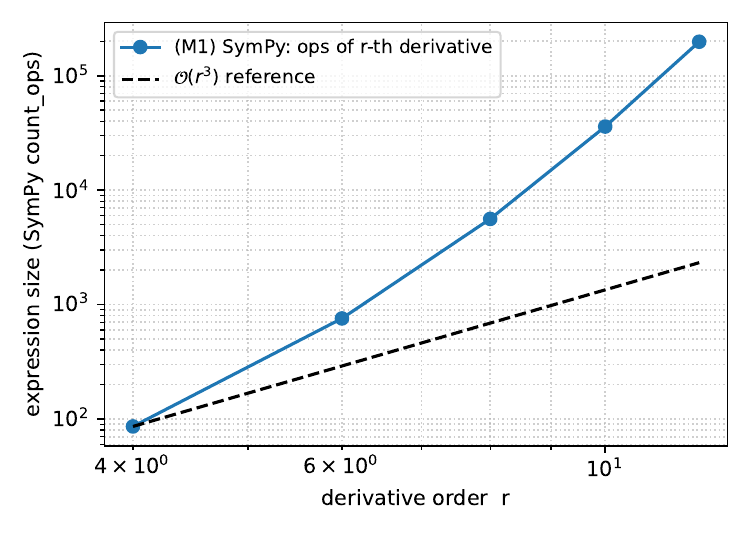}
\caption{Symbolic expression size of the $r$-th derivative built by~(M1), measured as \texttt{sympy.count\_ops}, versus the $O(r^{3})$ arithmetic-cost reference of~(M3).}
\label{fig:bench-size}
\end{figure}

\paragraph*{Results.} The full numerical results are gathered in Table~\ref{tab:bench} and plotted in Figs.~\ref{fig:bench-time} and~\ref{fig:bench-size}.

\begin{table}[t]
\centering
\caption{Benchmark of (M1) direct SymPy differentiation, (M2) Fa\`a di Bruno via Bell polynomials, (M3) formal-series composition (proposed), as a function of the requested derivative order~$r$. Times are in seconds (CPU). $\mathrm{ops}_{\mathrm{M1}}$ counts the AST nodes of the $r$-th symbolic derivative produced by (M1). The last column reports the max relative error $\max_{k<r} |h^{(k)}_{\mathrm{M2}}-h^{(k)}_{\mathrm{M3}}|/\max(|h^{(k)}_{\mathrm{M2}}|,|h^{(k)}_{\mathrm{M3}}|,1)$. Em-dashes indicate timeouts (20\,s for M1, 60\,s for M2).}
\label{tab:bench}
\input{Codes/data/benchmark_results.tex}
\end{table}

\paragraph*{Discussion.} Two walls are clearly visible.
\emph{First}, the symbolic method~(M1) grows super-polynomially: the size of the $r$-th derivative of $\cos g(s)$ contains, by virtue of Fa\`a di Bruno's formula, one term per integer partition of~$r$, weighted by products of derivatives of~$g$. With $g$ a generic polynomial of degree~$r$, none of these terms cancel, so SymPy's expression tree grows in lockstep with $p(r)$ and breaches the $20\,\mathrm{s}$ wall near $r\!=\!14$ on our hardware.
\emph{Second}, the Bell-polynomial route~(M2) decouples symbolic algebra from arithmetic but still has to enumerate partitions through the recursion~\eqref{eq:bell}; the un-memoised implementation \cite{Code} still completes $r\!=\!24$ in $7.7\,\mathrm{s}$ but exceeds its $60\,\mathrm{s}$ timeout by $r\!=\!32$. Even a memoised variant has $O(r^{3})$ arithmetic complexity but with a constant that is an order of magnitude larger than (M3) because each Bell evaluation re-reads $O(r)$ derivatives.
The proposed method~(M3) trades the partition-based combinatorics of \eqref{eq:bell} for two truncated polynomial compositions \eqref{eq:p123}, each of cost $O(r^{2})$ for the convolution and $O(r^{3})$ for the full composition. Empirically, the time of~(M3) follows the dashed cubic reference of Fig.~\ref{fig:bench-size} and remains under one millisecond up to $r=40$.

\paragraph*{Numerical accuracy.} In the regime where (M2) completes, the relative-error column of Table~\ref{tab:bench} remains at $\sim 10^{-15}$, i.e.\ machine precision in IEEE~754 double. The internal rescaling $f^{[k]}=f^{(k)}/k!$ used by~(M3) is what keeps the floating-point representation well conditioned: in the unscaled variables $h^{(k)}$ alone grows factorially and the absolute error scales with $k!$, but no precision is lost where it matters. This substantiates the ``without any loss of accuracy'' statement in the abstract.

\paragraph*{Impact on the trailer chain.} For the full pipeline of Section~II, each trailer~$k$ requires one block of cost~$O((n-k)^{3})$, so the total cost of~\eqref{eq:phi-factorised} is $O(n^{4})$ in the proposed method, versus the super-exponential cost of~(M1) and the $O(n\cdot p(n))$ partition-count cost of~(M2). Empirically, the full $n=20$ pipeline runs in well under a second on a laptop.

\section{Simulations}

\subsection{Animation speed}
\begin{figure}
\begin{center}
\begin{tabular}{cc}
\begin{minipage}{4cm}\begin{center}\includegraphics[width=1.2\textwidth]{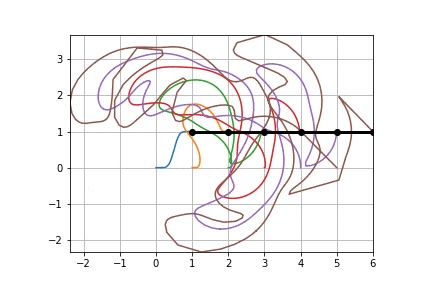}\end{center}\end{minipage}&
\begin{minipage}{4cm}\begin{center}\includegraphics[width=1.2\textwidth]{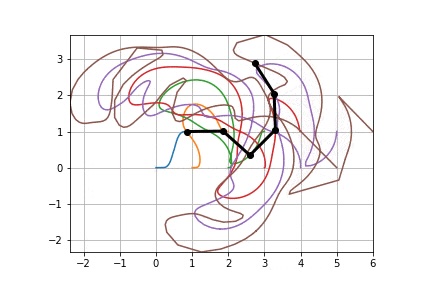}\end{center}\end{minipage}\\
\begin{minipage}{4cm}\begin{center}\includegraphics[width=1.2\textwidth]{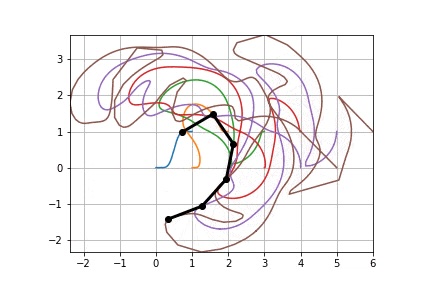}\end{center}\end{minipage}&
\begin{minipage}{4cm}\begin{center}\includegraphics[width=1.2\textwidth]{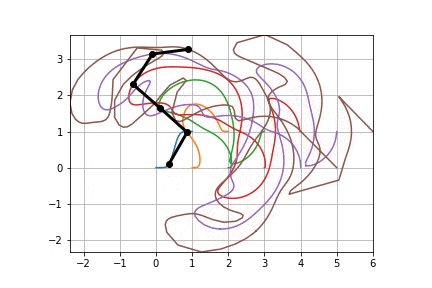}\end{center}\end{minipage}\\
\begin{minipage}{4cm}\begin{center}\includegraphics[width=1.2\textwidth]{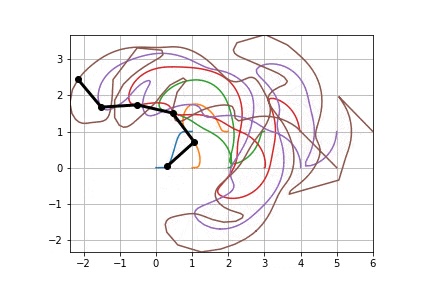}\end{center}\end{minipage}&
\begin{minipage}{4cm}\begin{center}\includegraphics[width=1.2\textwidth]{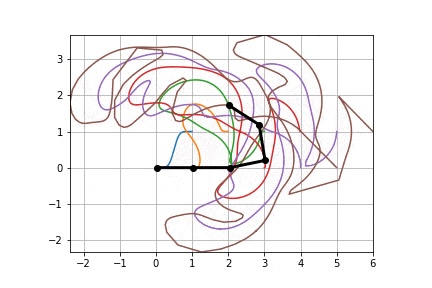}\end{center}\end{minipage}\\
\begin{minipage}{4cm}\begin{center}\includegraphics[width=1.2\textwidth]{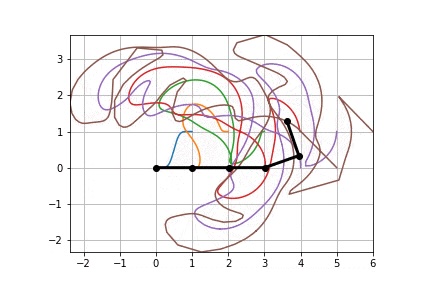}\end{center}\end{minipage}&
\begin{minipage}{4cm}\begin{center}\includegraphics[width=1.2\textwidth]{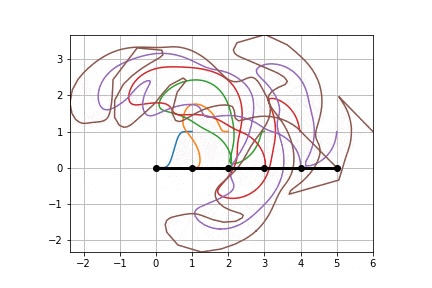}\end{center}\end{minipage}
\end{tabular}
\end{center}
\caption{Displacement and parking of a car with 5 trailers. The middle of each trailer's axle is marked by a black dot; each dot follows its own path (one colour per axle), and consecutive axles are linked by black segments. The initial configuration is top-left, the final one bottom-right; intermediate snapshots are arranged in row-major order. Initial and final configurations are close together and orderly, which is what produces the involute-looking trajectories needed to reach the goal.\label{simuFig}}
\end{figure}

To produce an animation of a car with many trailers, one chooses a trajectory $\boldsymbol{q}_1(s)$ for the last trailer's axle and uses the map $\varphi$ to compute all other trailer and car positions. Setting $s=t$ directly leads to a jerky animation because the head car's speed varies wildly (high derivatives) over time. We instead impose the head car's rear-axle speed (which can be computed from the derivatives of $\boldsymbol{q}_1(s(t))$) by integrating
\begin{equation}
\dot s(t) = \frac{1}{\sqrt{\bigl(x_{n-1}'(s(t))\bigr)^{2} + \bigl(y_{n-1}'(s(t))\bigr)^{2}}}.
\end{equation}
A simulation is shown in Fig.~\ref{simuFig} for $n=5$ trailers.

\section{Conclusion}

Differential flatness enables flexible reconfiguration of a car with trailers by exploiting the properties of flat outputs. Initial and final conditions define respectively the starting and ending configurations of the vehicle (with no singularities in the sense of~\cite{FLMR0} except at these boundary configurations). Feasible motions are produced by smooth interpolation between these states, typically with high-order polynomials, whose derivatives are then combined to compute the steering inputs.

This computation, performed naively, becomes intractable when the number of trailers~$n$ increases, due to the algebraic blow-up documented in Section~IV: direct symbolic differentiation hits a 20\,s wall at $r=14$ on our hardware, and Fa\`a di Bruno's formula via Bell polynomials hits a 60\,s wall beyond $r=24$. We have shown that an $O(n^{4})$ algorithm is obtained by combining four ingredients:
\begin{enumerate}
\item Suitable intermediate angular coordinates~\eqref{eq:angle-def}, which replace the offending square-root by the rational expression~\eqref{eq:alpha-prime};
\item Factorisation~\eqref{eq:phi-factorised} of the iterated map into single-trailer sub-maps, each operating on numerical arrays of derivatives rather than on symbolic expressions;
\item Convolutional, rescaled product and quotient rules~\eqref{eq:leibniz-rescaled}--\eqref{eq:recip} in place of the binomial-weighted Leibniz form;
\item Higher-order composition via truncated formal power series~\eqref{eq:p123} in place of the partition-based Fa\`a di Bruno formula.
\end{enumerate}
The resulting routine evaluates derivative orders up to $r=40$ in under one millisecond and matches the reference methods to machine precision where they succeed. The animation of Section~V is one illustration; an animation of a $20$-trailer manoeuvre and the full implementation code are available at~\cite{Code}.

\paragraph*{Future work.} Two extensions are immediate: (i)~replacing the heuristic time-reparametrisation of Section~V by a constraint-aware one that bounds the steering velocity and acceleration, and (ii)~using Algorithm~\ref{alg:main} inside an MPC outer loop so as to handle disturbances and modelling errors in real time. We expect the $O(n^{4})$ cost to remain well within real-time budgets for $n\le 20$.

\bibliography{papelard}
\end{document}

%% file: Codes/data/benchmark_results.tex
\begin{tabular}{rrrrrr}
\toprule
$r$ & $t_{\mathrm{M1}}$\,[s] & ops$_{\mathrm{M1}}$ & $t_{\mathrm{M2}}$\,[s] & $t_{\mathrm{M3}}$\,[s] & rel.\,err.\,$\mathrm{M2{-}M3}$ \\
\midrule
 4 & 1.40e-02 & 86 & 2.50e-05 & 4.30e-05 & 0 \\
 6 & 4.90e-02 & 756 & 4.20e-05 & 7.40e-05 & 4.06e-16 \\
 8 & 2.05e-01 & 5\,601 & 1.55e-04 & 1.15e-04 & 4.94e-16 \\
10 & 1.05e+00 & 35\,917 & 5.34e-04 & 1.83e-04 & 9.72e-16 \\
12 & 5.06e+00 & 197\,978 & 1.88e-03 & 2.77e-04 & 1.59e-15 \\
14 & \textemdash & \textemdash & 7.43e-03 & 4.03e-04 & 1.40e-15 \\
16 & \textemdash & \textemdash & 3.06e-02 & 5.08e-04 & 1.52e-15 \\
20 & \textemdash & \textemdash & 4.83e-01 & 8.95e-04 & 7.93e-16 \\
24 & \textemdash & \textemdash & 7.74e+00 & 1.35e-03 & 9.05e-16 \\
32 & \textemdash & \textemdash & \textemdash & 2.69e-03 & \textemdash \\
40 & \textemdash & \textemdash & \textemdash & 4.65e-03 & \textemdash \\
\bottomrule
\end{tabular}